\documentclass[reqno]{amsart}
\usepackage{graphicx}
\usepackage{epsfig}
\usepackage[latin1]{inputenc}
\usepackage{amsmath}
\usepackage{amsfonts}
\usepackage{amssymb}
\usepackage{amsthm}
\usepackage{amsmath, amscd}
\usepackage{verbatim}
\usepackage{subfigure}

\usepackage{verbatim}
	
\newtheorem*{Theorem}{Theorem}

\def\rk{\mathrm{rk}}
\def\Z{\mathbb{Z}}

\def\P{\mathcal{P}}

\title{A Note on Cabling and L-space surgeries}

\subjclass[2009]{}
\author{Jennifer Hom}

\address{Department of Mathematics, University of Pennsylvania, Philadelphia, PA 19104
\newline\indent{\tt jenhom@math.upenn.edu}}

\numberwithin{equation}{section}

\begin{document}

\begin{abstract}
We prove that the $(p,q)$-cable of a knot $K \subset S^3$ admits a positive $L$-space surgery  if and only if $K$ admits a positive $L$-space surgery and $q/p\geq 2g(K)-1$, where $g(K)$ is the Seifert genus of $K$. The ``if" direction is due to Hedden \cite{HeddencablingII}.
\end{abstract}

\maketitle

\section{Introduction}

In \cite{OS3manifolds}, Ozsv\'{a}th and Szab\'{o} introduced a powerful tool for studying closed $3$-manifolds, Heegaard Floer homology, and later equipped this invariant with a filtration (independently developed by Rasmussen in \cite{R}) that defined an invariant for a knot in the $3$-manifold. The relationship between the knot invariant and the Heegaard Floer homology of the $3$-manifold obtained by Dehn surgery on that knot has been well studied (c.f. \cite{OSlens}, \cite{OSinteger}, \cite{OSrational}), and can also be considered from the perspective of bordered Heegaard Floer homology \cite{LOT}.

In this note, we restrict our consideration to the simplest ``hat" version of the theory, assuming that the reader is familiar with the finitely generated abelian groups
$$
\widehat{HF}(Y) \qquad\textup{and}\qquad \widehat{HFK}(Y,K)
$$
associated with a $3$-manifold $Y$ and a null-homologous knot $K\subset Y$ (see \cite{OSknots}). We will write simply $\widehat{HFK}(K)$ when it is clear that the ambient $3$-manifold is $S^3$. For the present purposes, we do not need to concern ourselves with the gradings on these groups.
We focus our attention on a class of $3$-manifolds with particularly simple Heegaard Floer homology. For a rational homology sphere $Y$, we have that 
$$\mathrm{rk} \ \widehat{HF}(Y) \geq | H_1 (Y, \Z)|.$$
An \emph{$L$-space} is a rational homology sphere $Y$ for which the above bound is sharp. The name comes from the fact that lens spaces are $L$-spaces, which can be seen by examining the Heegaard Floer complex associated to a standard genus $1$ Heegaard decomposition of a lens space.

We call a knot $K \subset S^3$ an \emph{$L$-space knot} if there exists $n\in \Z$, $n >0$, such that $n$ surgery on $K$ yields an $L$-space. We will denote the resulting $3$-manifold by $S_n^3(K)$. Torus knots are a convenient source of $L$-space knots, since $pq \pm 1$ surgery on the $(p,q)$-torus knot yields a lens space. It was proved in \cite{OSlens} that if a knot $K$ is an $L$-space knot, then the knot Floer complex associated to $K$ has a particularly simple form that can be deduced from the Alexander polynomial of $K$, $\Delta_K (t)$. Thus, knowing that a knot $K$ admits a lens space (or $L$-space) surgery yields a remarkable amount of information about the Heegaard Floer invariants associated to both the knot $K$, and manifolds obtained by Dehn surgery on $K$. In particular, results of \cite{OSlens} and \cite{OSrational} allow one to compute the Heegaard Floer invariants of any Dehn surgery on an $L$-space knot $K$ from the Alexander polynomial of $K$.

Recall that the $(p,q)$-cable of a knot $K$, denoted $K_{p,q}$, is the satellite knot with pattern the $(p,q)$-torus knot. More precisely, we can construct $K_{p,q}$ by equipping the boundary of a tubular neighborhood of $K$ with the $(p,q)$-torus knot, where the knot traverses the longitudinal direction $p$ times and the meridional direction $q$ times. We will assume throughout that $p>1$. (This assumption does not cause any loss of generality, since $K_{-p,-q}=-K_{p,q}$, where $-K_{p,q}$ denotes $K_{p,q}$ with the opposite orientation, and since $K_{1,q}=K$.) Let $g(K)$ denote the Seifert genus of $K$. In Theorem 1.10 of \cite{HeddencablingII}, Hedden proves that if $K$ is an $L$-space knot and $q/p \geq 2g(K)-1$, then $K_{p,q}$ is an $L$-space knot. The goal of this note is to prove the converse:

\begin{Theorem}
\label{main}
The $(p,q)$-cable of a knot $K \subset S^3$ is an $L$-space knot if and only if $K$ is an $L$-space knot and $q/p\geq 2g(K)-1$.
\end{Theorem}
\noindent We prove our theorem by methods similar to those used in \cite{HeddencablingII}. An interesting question to consider is whether there are other satellite constructions that also yield $L$-space knots.
\newline \newline
\noindent \textbf{Acknowledgements.} I would like to thank Paul Melvin for his patience and encouragement, and Matt Hedden for helpful email correspondence.

\section{Proof of Theorem}

An $L$-space $Y$ can be thought of as rational homology sphere with the ``smallest" possible Heegaard Floer invariants, i.e. $\rk \ \widehat{HF}(Y)=|H_1(Y, \Z)|$. In a similar spirit, an $L$-space knot $K$ can be thought of as a knot with the ``smallest" possible knot Floer invariants. For example, since
$$\Delta_K(t)=\sum_{m,s} (-1)^m \ \rk \ \widehat{HFK}_m(K,s)t^s,$$
we see immediately that the total rank of $\widehat{HFK}(K)$ is bounded below by the sum of the absolute value of the coefficients of the Alexander polynomial of $K$, $\Delta_K(t)$. A necessary, but not sufficient, condition for a knot $K$ to be an $L$-space knot is for this bound to be sharp; see \cite{OSlens} for the complete statement. The spirit of our proof is that when either $K$ is not an $L$-space knot, or $q/p<2g(K)-1$, the knot Floer invariants of $K_{p,q}$ are not ``small" enough for $K_{p,q}$ to be an $L$-space knot. We will determine this by looking at the rank of $\widehat{HF}(S_{pq}^3(K_{p,q}))$.

Recall that $\tau(K)$ is the integer-valued concordance invariant defined by Ozsv\'{a}th and Szab\'{o} in \cite{OS4ball}. Let $\P$ denote the set of all knots $K$ for which $g(K)=\tau(K)$. We begin by assembling the following collection of facts. 

\begin{enumerate}
\item \label{LP} If $K$ is an $L$-space knot, then $K \in \P$ \cite{OSlens}.
\item \label{prop95} Let $$s_K=\sum_{s \in \Z}(\rk \ H_*(\hat{A}_s(K))-1),$$ where $\hat{A}_s(K)$ is the subcomplex of $CFK^{\infty}(K)$ defined in \cite{OSrational}. We may think of $CFK^{\infty}(K)$ as generated over $\Z[U, U^{-1}]$ by $\widehat{HFK}(K)$, in which case $\rk \ \hat{A}_s(K)=\rk \ \widehat {HFK}(K)$ for all $s$. Recall that $\rk \ \widehat {HFK}(K)$ is always odd, since the graded Euler characteristic of $\widehat{HFK}(K)$ is the Alexader polynomial of $K$. Therefore, $\rk \ H_*(\hat{A}_s(K))$ is odd, hence greater than or equal to $1$, and so $s_K$ is always non-negative. Let $$t_{K}^{a/b}=2\ \max(0,(2g(K)-1)b-a),$$ for a pair of relatively prime integers $a$ and $b$, $b>0$. Notice that $$t_{K}^{a/b}=0 \textup{ if and only if } a/b \geq 2g(K)-1.$$
For $K \in \P$ and $a$, $b$ as above,
$$\rk \ \widehat{HF}(S_{a/b}^3(K))=a+bs_K+t_{K}^{a/b}.$$
This is a special case of Proposition 9.5 of \cite{OSrational}. In particular, the term $\nu(K)$ appearing in Proposition 9.5 is bounded below by $\tau(K)$ and above by $g(K)$, so $K \in \P$ implies $\nu(K)=g(K)$. We notice that
$$K  \textup{ admits a positive $L$-space surgery if and only if } s_K=0.$$ If $s_K=0$, then $p$ surgery on $K$ yields an $L$-space, for any integer $p \geq 2g(K)-1$, and if $K$ is an $L$-space knot, then there exists some integer $p>0$ such that $p$ surgery on $K$ is an $L$-space, in which case $s_K$, which is always non-negative, must be $0$.
\item \label{t} Recall our convention that $p$, $q$ are relatively prime integers, with $p>1$. If $K_{p,q}\in \P$, then $K \in \P$, and if $K \in \P$, then $\tau(K_{p,q})=p\tau(K)+\tfrac{(p-1)(q-1)}{2}$.
These facts are Corollaries 4 and 3, respectively, in \cite{VanCott}. Therefore, if $K_{p,q}\in \P$, we have
\begin{align*}
(2g(K)-1)p-q&= (2\tau(K)-1)p-q \\
&= 2\big(p\tau(K)+\tfrac{(p-1)(q-1)}{2}\big)-1-pq \\
&= 2\tau(K_{p,q})-1-pq\\
&= 2g(K_{p,q})-1-pq,
\end{align*}
or equivalently, 
$$\textup{if } K_{p,q} \in \P, \textup{ then }t_{K}^{q/p}=t_{K_{p,q}}^{pq}.$$
\item \label{pq} $pq$ surgery on $K_{p,q}$ is the manifold $L(p,q)\# S_{q/p}^3(K)$. See \cite{HeddencablingII} for a nice proof of this classical fact about Dehn surgery. We also have from \cite{OS3manifolds} that
$$\rk \ \widehat{HF}(Y_1 \# Y_2)= \rk \ \widehat{HF}(Y_1)\cdot \rk \ \widehat{HF}(Y_2).$$
Then
\begin{align*}
\rk \ \widehat{HF}(S_{pq}^3(K_{p,q}))&= \rk \ \widehat{HF}(L(p,q)) \cdot \rk \ \widehat{HF}(S_{q/p}^3(K))  \\
&= p \cdot \rk \ \widehat{HF}(S_{q/p}^3(K)).
\end{align*}
\end{enumerate}

\noindent With these facts in place, we are ready to prove the theorem. Assume $K_{p,q}$ is an $L$-space knot. Then by (\ref{LP}) and (\ref{t}), $K_{p,q} \in \P$ and $t_{K_{p,q}}^{pq}=t_{K}^{q/p}$, and by (\ref{prop95}),
$$\rk \ \widehat{HF}(S_{pq}^3(K_{p,q}))=pq+s_{K_{p,q}}+t_{K_{p,q}}^{pq}  \qquad\textup{and}\qquad \rk \ \widehat{HF}(S_{q/p}^3(K))=q+ps_K+t_{K}^{q/p}.$$
Then by (\ref{pq}), $\rk \ \widehat{HF}(S_{pq}^3(K_{p,q}))=p\cdot \rk \ \widehat{HF}(S_{q/p}^3(K))$, and $s_{K_{p,q}}=0$, since $K_{p,q}$ is an $L$-space knot. So we find that
$$p^2s_K+(p-1)t_{K}^{q/p}=0.$$
Therefore, since $p>1$, we have that $s_K$ and $t_{K}^{q/p}$ must both be zero, or equivalently, $K$ is an $L$-space knot and $q/p \geq 2g(K)-1$. This completes the proof of the theorem.

\end{document}